\documentclass[10pt, oneside]{amsart}
\usepackage{amsmath,amsthm,amsfonts,amscd, 
amssymb, epsf}
\usepackage{url}
\usepackage{geometry}
\makeatletter

\theoremstyle{plain}
\newtheorem{thm}{Theorem}[section]
  \theoremstyle{plain}
  \newtheorem{lem}[thm]{Lemma}
   \newtheorem{cor}[thm]{Corollary}
   
  \newtheorem{as}[thm]{Assumption}
  \newtheorem*{thm*}{Theorem}
  \newtheorem*{cor*}{Corollary}
\numberwithin{equation}{section}

\usepackage{palatino}
\usepackage{amsthm}
\usepackage{amsfonts}
\usepackage{amscd}
\usepackage{epsf}
\makeatletter

\newcommand{\field}[1]{\ensuremath{\mathbb
{#1}}}

\newcommand{\CC}{\field{C}}
\newcommand{\df}{\equiv}
\newcommand{\NN}{\field{N}}

\newcommand{\HHHH}{\field{H}}
\newcommand{\hh}{\field{H} \,}

\DeclareMathOperator{\HH}{\HHHH^3}

\newcommand{\RR}{\field{R}}

\newcommand{\ZZ}{\field{Z}}

\newcommand{\F}{\mathcal{F}}
\newcommand{\E}{\mathcal{E}}
\newcommand{\scz}{\mathcal{S}}

\newcommand{\K}{\mathcal{K}}
\newcommand{\D}{\mathcal{D}}
\newcommand{\hil}{\mathcal{H}}

\newcommand{\M}{\mathcal{M}}

\DeclareMathOperator{\rep}{Rep}
\DeclareMathOperator{\PSL}{PSL}

\DeclareMathOperator{\vol}{vol}

\DeclareMathOperator{\R}{Re}
\DeclareMathOperator{\I}{Im}
\DeclareMathOperator{\pc}{PSL(2,\CC)}
\DeclareMathOperator{\GL}{GL}

\DeclareMathOperator{\tr}{tr}

\DeclareMathOperator{\en}{\mathnormal{\mathcal
{E}(T)}}
\DeclareMathOperator{\oen}{\mathnormal{\left|
\mathcal{E}(T) \right|}}

\DeclareMathOperator{\ren}{\mathnormal
{\mathcal{E}(R)}}

\DeclareMathOperator{\hs}{\mathnormal{\hil
(\Gamma,\chi)}}
\DeclareMathOperator{\lp}{\mathnormal
{\varDelta}}

\DeclareMathOperator{\com}{\textbf{Comm}}
\newcommand{\beq}{\begin{equation}}
\newcommand{\eeq}{\end{equation}}

\newcommand{\C}{\mathcal{C}}

\newcommand{\ra}{\rightarrow}

\makeatother

\begin{document}
\noindent
\author{Joshua S. Friedman}
\address{ Department of Mathematics and 
Sciences, United States Merchant Marine 
Academy, 300 Steamboat Road, Kings Point, NY  
11024}
\email{CrownEagle@gmail.com}
\email{friedmanj@usmma.edu}
\email{joshua@math.sunysb.edu}

\title[Modular Correspondences]{The Selberg 
Trace Formula for Hecke operators on cocompact Kleinian groups}

\begin{abstract}
We compute the Selberg trace 
formula for Hecke operators (also called the trace formula for modular correspondences) in the 
context of cocompact Kleinian groups with 
finite-dimentional unitary representations. We 
give some applications to the distribution of Hecke eigenvalues, and give an analogue of Huber's theorem.
\end{abstract}

\maketitle

\thispagestyle{empty}

\section{Introduction}
\subsection{Motivation} The Selberg trace formula has been well  studied for cofinite Kleinian and Fuchsian groups. In the Fuchsian case: \cite{Selberg1, Roelcke, VenKalFad,  Hejhal1, Hejhal2, Venkov, Buser, Iwaniec}; in the Kleinian case: \cite{Elstrodt, Friedman1, Friedman2}. For the most part (excepting \cite{Hejhal1,Hejhal2}) the references above consider the \emph{standard} Selberg trace formula which entails two objects: a discrete group $\Gamma,$ and possibly a finite-dimensional unitary representation $\chi.$ A third object to consider is a Hecke operator. 

The Selberg trace formula (and the underlying Selberg spectral theory) for Hecke operators (also known as the Selberg trace formula for \emph{modular correspondences})\footnote{Some authors seem to refer to the subject as \emph{The Selberg trace formula for modular correspondences} and some as \emph{The Selberg trace formula for Hecke operators.} For our purposes these are interchangeable. To be perfectly precise, we are not computing the trace of a Hecke operator, but rather the trace of $\M h(\lp),$ where $\M$ is a Hecke operator, $\lp$ is the Laplacacian, and $h$ is a holomorphic function with a certain decay rate.  } is far from a trivial extension of the standard Selberg trace formula. By considering the abundance of Hecke operators for the Modular group $\Gamma = \PSL(2,\ZZ),$  Lindenstrauss and Venkatesh (\cite{Linden}, \cite[Chap. 4]{Goldfeld}) gave a new proof\footnote{The old proof is Selberg's original method: The determinant of the scattering matrix is shown to have order that is sufficiently small.  } for  the existence of infinitely many even Maa\ss\,  forms\footnote{Their proof works in much more generality.} in $L^{2}(\hh \setminus \Gamma).$    

Str\"{o}mbergsson (\cite{Strom2}) made extensive use of the trace formula for Hecke operators to solve a problem  in the \emph{well known spectral correspondence for quaternion groups}---the Jacquet-Langlands correspondence (\cite{Hejhal4}). More specifically, Str\"{o}mbergsson completely determined the image of the Jacquet-Langlands correspondence.\footnote{The bijective correspondence between the nontrivial automorphic forms on the multiplicative group of a division quaternion algebra and certain cusp forms on $\GL(2).$} In fact, in \cite{Strom1}, he extended the Selberg trace formula for modular correspondences from the cocompact Fuchsian case (\cite{Hejhal1}) to the cofinite Fuchsian case.

The possibility for a  trace formula for Hecke operators  started with Selberg \cite[pp. 444--446, 460--462, 504--505]{Selberg2}. In \cite{Hejhal1} Hejhal gave an explicit trace formula for modular correspondences for  cocompact Fuchsian groups with finite-dimensional unitary representations and integer weight $k \geq 0$ forms. Akiyama and Tanigawa \cite{Akiyama} proved the Selberg trace formula for modular correspondences for cofinite Fuchsian groups and gave explicit evaluations of the trace formula for the case of the congruence group $\Gamma_{0}(p)$ ($p$ prime). Hoffmann (\cite{Hoffmann1, Hoffmann2}) approached the trace formula for Hecke operator, from the representation-theoretic point of view.

Str\"{o}mbergsson (\cite{Strom1}) gave a very detailed derivation of the Selberg trace formula for modular correspondences in the context of weight-zero forms (functions) and trivial unitary representations on cofinite Fuchsian groups. He also gave explicit formulas for the case of $\Gamma = \Gamma_{0}(N)$ ($N$ is square-free).

In this work, we derive the Selberg trace formula for Hecke operators for cocompact Kleinian groups with finite-dimensional unitary representations. In order for the formula to apply, the cocompact group $\Gamma$ needs to have its commensurator $\com(\Gamma)$ (in $\pc$) strictly larger than itself. For each $\alpha \in \com(\Gamma) \setminus \Gamma,$ we obtain a Hecke operator $\M.$ When $\Gamma$ is \emph{arithmetic,} $[\com(\Gamma):\Gamma] = \infty;$ in fact $\com(\Gamma)$ is a dense subgroup of $\pc$ (\cite[p. 271, Exer. 8.4.5]{Reid}).

\subsection{Main results}
In this section we state our main results, and  briefly describe the notation and preliminary material needed to state them (for more details see \S\ref{secPre}).

Let $\Gamma < \pc $ be a cocompact Kleinian 
group acting on hyperbolic three-space $\HH.$ 
Let $V$ be a finite-dimensional complex inner 
product space with inner-product $\left<~,~
\right>_{V},$ and let $\rep(\Gamma,V)$ denote 
the set of finite-dimensional unitary 
representations of $\Gamma$ in $V.$ Let $\F 
\subset \Gamma$ be a (compact) fundamental 
domain for the action of $\Gamma$ in $\HH.$ 

Let $\chi \in \rep(\Gamma,V).$ The Hilbert 
space of $\chi-$\emph{automorphic} functions 
is the set of measurable functions    
\begin{multline*}
\hs \df  \{ f: \HH \ra V ~|~ f(\gamma P) = 
\chi(\gamma) f(P)~\forall \gamma \in \Gamma,  
P \in \HH, \\ $ and $ \left<f,f \right> \df 
\int_{\F} \left<f(P),f(P)\right>_V\,dv(P) < 
\infty \}. 
\end{multline*}
Finally, let   $ \lp = \lp(\Gamma,\chi) $ be 
the corresponding positive self-adjoint 
Laplace-Beltrami operator on $\hs.$

The \emph{commensurator} subgroup 
$\com(\Gamma)$ is the set of all $\alpha \in 
\pc $ such that both $[\Gamma : \Gamma \cap 
\alpha^{-1}\Gamma \alpha] < \infty$  and $
[\alpha^{-1}\Gamma \alpha : \Gamma \cap 
\alpha^{-1}\Gamma \alpha] < \infty. $ 

Throughout this paper we will consider elements $\alpha \in \com(\Gamma)$ that satisfy:
\begin{enumerate}
\item the element $\alpha \in \com(\Gamma) \setminus \Gamma;$
\item the unitary representation $\chi \in \rep
(\Gamma,V)$ has a single-valued extension from 
$\Gamma$ to the set $\Gamma \alpha \Gamma,$ 
and satisfies the following properties:
\begin{enumerate}
\item $\chi(g_{1} \alpha g_{2}) = \chi(g_{1}) \chi(\alpha) \chi(g_{2})$ for $g_{1},g_{2}\in \Gamma;$
\item $\chi(\alpha)$ is an invertible linear map on the inner product space $V.$ 
\end{enumerate}
\end{enumerate}

Associated to the element $\alpha$ is the \emph{Hecke operator} $\M:\hs \mapsto \hs.$ The Hecke operator commutes with $\lp.$

Let $\D$ be an indexing set for the set of eigenvalues of $\lp.$ For each eigenvalue  $\lambda = \lambda_{m} ~ (m \in \D),$ let $A_{\lambda} \df A(\lambda, \Gamma, \chi)$ denote  the subspace of $\hs$ spanned by $\{ \lambda_{m}~|~ \lambda_{m} = \lambda   \}.$ The subspaces $A_{\lambda}$ are invariant under $\M.$ Let $\lambda_{m_{1}}\dots \lambda_{m_{k_{\lambda}}} $ generate $A_{\lambda}$ ($\lambda$ has multiplicity $k_{\lambda}$). Since $A_{\lambda}$ is invariant under $\M,$ the eigenvalues of $\left. \M \right|_{A_{\lambda}}$ can be listed as $\omega_{m_{1}}\dots \omega_{m_{k_{\lambda}}}.$ 
Decomposing $\hs$ into invariant subspaces $A_{\lambda},$ for each $m \in \D,$ we pair $\lambda_{m}$ and $\omega_{m}.$ 

We next state our main result, the Selberg trace formula for Hecke operators. The various notations are described very briefly following the theorem; for more details see \S\ref{secSelTra}.

Our main result is the following:
\begin{thm*}{(Selberg trace formula)} 
Let $\Gamma $ be a cocompact Kleinian group. Let \emph{$\alpha \in\com(\Gamma) \setminus \Gamma,$} and let  $\chi$ be a representation satisfying Assumption~\ref{asAlpha}. Let  $h$ be a holomorphic function on  $ \{ s \in \CC \, | \, |\I(s)| < 2+ \delta \}$ for some $\delta > 0,$ satisfying $ h(1+z^2) = O( 1+|z|^2)^{3/2 - \epsilon}) $ as  $|z| \ra \infty,$ and let
$$ g(x) = \frac{1}{2\pi} \int_{\RR} h(1+t^2)e^{-itx}\,dt. $$  Then
\begin{multline*}
\sum_{m \in \D} h(\lambda_m) \omega_{m} = \sum_{ \{R \}_{\text{\emph{ell}}}} \frac{\tr_{V}\!\left( \chi(R^{-1})^{*} \right)   g(0) \log N(T_0) }{  |\E(R)| \, |(\tr(R))^{2}-4|}  + 
\sum_{\{ T \}_\text{\emph{lox}} } \frac{\tr_{V}\!\left( \chi(T^{-1})^{*} \right) g(\log N(T))}{\oen |a(T)-a(T)^{-1}|^{2}}\log N(T_{0})
\end{multline*}
\end{thm*}
Here, $ \{ \lambda _{m} \}_ {m \in \D } $ are the eigenvalues of $\lp$ counted with multiplicity, and $\omega_{m}$ are the eigenvalues of $\M$ (the Hecke operator associated to $\alpha$) with the convention of Equation~\ref{eqTraceMod}.  The summation with respect to $\{R \}_\text{ell} $ extends over the finitely many $\Gamma-$conjugacy classes of  elliptic elements  $R \in \Gamma \alpha^{-1} \Gamma,$ and for such a class, $N(T_0)$ is the minimal norm of a hyperbolic or loxodromic element of the centralizer $\mathcal{C}(R) \subset \Gamma;$ $\ren$ is the maximal finite subgroup contained in $ \mathcal{C}(R)$ (see Lemma~\ref{lemEllip} for more details).  The summation with respect to $\{ T \}_\text{lox}$ extends over the  $\Gamma-$conjugacy classes of hyperbolic or loxodromic elements of $\Gamma \alpha^{-1} \Gamma,$ $T_0$ denotes a primitive hyperbolic or loxodromic element of minimal norm in\footnote{Note that there is no clear relationship between $T \in \Gamma \alpha^{-1} \Gamma$ and $T_{0} \in \Gamma.$} $\mathcal{C}(T) \subset \Gamma;$  the element $T$ is conjugate in $\pc$ to the transformation  described by the diagonal matrix with diagonal entries $a(T), a(T)^{-1}$ with $|a(T)| > 1, $ and\footnote{Please note that there is a typographical error in the loxodromic and non cuspidal elliptic terms in \cite{Elstrodt} Theorem 6.5.1; both terms are missing a factor of $\frac{1}{4 \pi}. $} $N(T) =   |a(T)|^2; $ $\en$ is the finite cyclic elliptic subgroup of $\mathcal{C}(T)$  (see Lemma~\ref{lemLoxo} for more details).  The sum over elliptic elements is finite, and all other sums converge absolutely.

Please keep in mind that the conjugacy classes $\{R\}_{\text{ell}}$ and $\{T\}_{\text{lox}},$  the operator $\M,$ and the eigenvalues $\omega_{m}$ all depend on $\alpha \in \com(\Gamma) \setminus \Gamma.$ 

Our first application of the trace formula is to study the distribution of the Hecke eigenvalues.

Define the \emph{elliptic number of $\Gamma$ with respect to $\alpha,$} $E_{\alpha},$  by the finite sum 
$$
E_{\Gamma}^{\alpha} \df \sum_{ \{R \}_{\text{{ell}}}} \frac{ \log N(T_0) }{  |\E(R)| \, |(\tr(R))^{2}-4|}. 
$$
Here the sum is over \emph{elliptic conjugacy classes} of the set $\Gamma \alpha^{-1} \Gamma,$ where conjugacy is defined with respect to the group $\Gamma.$ The notations above are defined in \S\ref{secEllCas}.  

We have 
\begin{thm*} 
Let $\Gamma$ be a cocompact Kleinian group with $\alpha \in \com(\Gamma) \setminus \Gamma.$  Then
$$  \sum_{m \in \D} \omega_{m} e^{-\lambda_{m}t} = \frac{E_{\Gamma}^{\alpha}}{\sqrt{4\pi t}} + O(\sqrt{t}) \quad \text{as $t \ra 0^{+}.$} $$
Here $E_{\Gamma}^{\alpha}$ is the elliptic number of $\Gamma$ with respect to $\alpha,$ $\{ \lambda_{m} \}_{m \in \D } $ are the eigenvalues of $\lp$ counted with multiplicity, and $\omega_{m}$ are the eigenvalues of $\M$ (the Hecke operator associated to $\alpha$) with the convention of Equation~\ref{eqTraceMod}.
\end{thm*}

If the set $\Gamma \alpha^{-1} \Gamma$ contains no elliptic elements we have
\begin{thm*} 
Let $\Gamma$ be a cocompact Kleinian group with $\alpha \in \com(\Gamma) \setminus \Gamma.$ Suppose that the set $\Gamma \alpha^{-1} \Gamma$ contains no elliptic elements. Then
$$  \sum_{m \in \D} \omega_{m} e^{-\lambda_{m}t} = 
 O\left(t^{-1/2}\exp(-c/t)\right)  
 \quad \text{as $t \ra 0^{+}.$} $$
Hence
$$
\lim_{t \ra 0^{+}} \sum_{m \in \D} \omega_{m} e^{-\lambda_{m}t}   = 0.$$
Here $\{ \lambda_{m} \}_{m \in \D } $ are the eigenvalues of $\lp$ counted with multiplicity, and $\omega_{m}$ are the eigenvalues of $\M$ (the Hecke operator associated to $\alpha$) with the convention of Equation~\ref{eqTraceMod}.
\end{thm*}

Finally, we give an analogue of Huber's theorem.

Following \cite{Elstrodt}, we define the length spectrum of $\Gamma \alpha^{-1} \Gamma.$ For loxodromic $T_{j} \in \Gamma \alpha^{-1} \Gamma$ set $\mu_{j} = \log N(T_{j}).$  The \emph{length spectrum} of $\Gamma \alpha^{-1} \Gamma$ is defined to be    
$$ \mathcal{L}_{\Gamma}^{\alpha} \df \left( \mu_{j}, \sum_{ \substack{ \{ T \}_\text{\emph{lox}} \\ \log N(T)=\mu_{j} } }  \frac{\log N(T_{0}) }{\oen |a(T)-a(T)^{-1}|^{2}}  \right)_{j \geq 1.} $$ 

In the two-dimensional case the length spectrum simply comprises the lengths of closed geodesics. Here we really need the \emph{complex lengths} $a(T)$ and the order of the elliptic, finite subgroup of $\C(T),$ $\oen.$  

We define the eigenvalue spectra of $\lp$ and $\M$ ($\M$ is defined from $\alpha \in  \com(\Gamma) \setminus \Gamma$) by
$$\mathcal{S}_{\Gamma}^{\alpha} = \left(\lambda_{j}, \omega(\lambda_{j}) \right)_{j \in \D^{*}}.$$
Here the symbol $D^{*}$ means that we do not count with multiplicity, and  $\omega(\lambda)$ is the trace of $\M$ on the invariant subspace generated by all eigenfunctions (of $\lp$) with eigenvalue $\lambda.$ 

\begin{thm*} Let $\Gamma, \Gamma'$ be cocompact Kleinian groups with $\alpha \in \com(\Gamma) \setminus \Gamma,~\alpha' \in \com(\Gamma') \setminus \Gamma'. $ Then the following hold:

\begin{enumerate}
\item Suppose that $\mathcal{S}_{\Gamma}^{\alpha}$ and $\mathcal{S}_{\Gamma'}^{\alpha'}$ agree up to at most finitely many terms. Then $$E_{\Gamma}^{\alpha} = E_{\Gamma'}^{\alpha'}$$ 
$$\mathcal{S}_{\Gamma}^{\alpha} = \mathcal{S}_{\Gamma'}^{\alpha'}$$  
$$\mathcal{L}_{\Gamma}^{\alpha} = \mathcal{L}_{\Gamma'}^{\alpha'}.$$
\item Suppose that $\mathcal{L}_{\Gamma}^{\alpha}$ and $\mathcal{L}_{\Gamma'}^{\alpha'}$ agree up to at most finitely many terms. Then $$E_{\Gamma}^{\alpha} = E_{\Gamma'}^{\alpha'}$$ 
$$\mathcal{S}_{\Gamma}^{\alpha} = \mathcal{S}_{\Gamma'}^{\alpha'}$$  
$$\mathcal{L}_{\Gamma}^{\alpha} = \mathcal{L}_{\Gamma'}^{\alpha'}.$$
\end{enumerate}
\end{thm*}

By applying the above Theorem to the case of a \emph{single} fixed group with two different elements $\alpha, \alpha' \in \com(\Gamma) \setminus \Gamma, $ we obtain:

\begin{cor*}
Let $\Gamma$ be a cocompact Kleinian groups with $\alpha, \alpha' \in \com(\Gamma) \setminus \Gamma. $ Let $\omega(\lambda_{m}) = \omega'(\lambda_{m})$ for all but at most finitely many $m \in \D^{*}.$ Then $\omega(\lambda_{m}) = \omega'(\lambda_{m})$ for all $m \in \D^{*}.$ 
\end{cor*}

\section{Preliminaries} \label{secPre}
In this section we review the basic notions 
needed to evaluate the Selberg trace formula 
for modular correspondences. See \cite[Chapter 
5]{Hejhal1} for the analogous Fuchsian case\footnote{There is much overlap between the 
Fuchsian case and our case, the Kleinian case. 
We will often cite the Fuchsian case, as a 
reference when proofs in the two cases are 
identical.}.

\subsection{Cocompact Kleinian groups} 
Let $\Gamma < \pc $ be a cocompact Kleinian 
group acting on hyperbolic three-space $\HH.$ 
Let $V$ be a finite-dimensional complex inner 
product space with inner-product $\left<~,~
\right>_{V},$ and let $\rep(\Gamma,V)$ denote 
the set of finite-dimensional unitary 
representations of $\Gamma$ in $V.$ Let $\F 
\subset \Gamma$ be a (compact) fundamental 
domain for the action of $\Gamma$ in $\HH.$ 

Let $\chi \in \rep(\Gamma,V).$ The Hilbert 
space of $\chi-$\emph{automorphic} functions 
is the set of measurable functions    
\begin{multline*}
\hs \df  \{ f: \HH \ra V ~|~ f(\gamma P) = 
\chi(\gamma) f(P)~\forall \gamma \in \Gamma,  
P \in \HH, \\ $ and $ \left<f,f \right> \df 
\int_{\F} \left<f(P),f(P)\right>_V\,dv(P) < 
\infty \}. 
\end{multline*}
Finally, let   $ \lp = \lp(\Gamma,\chi) $ be 
the corresponding positive self-adjoint 
Laplace-Beltrami operator on $\hs.$

\subsection{Hecke operators}
In this section we define the Hecke operators (or modular correspondences). For more details see \cite[Chap. 5]{Hejhal1} for the cocompact Fuchsian case, \cite{Strom1} for the cofinite Fuchsian case.

Let $\Gamma$ be a cocompact Kleinian group. 
The \emph{commensurator} subgroup 
$\com(\Gamma)$ is the set of all $\alpha \in 
\pc $ such that both $[\Gamma : \Gamma \cap 
\alpha^{-1}\Gamma \alpha] < \infty$  and $
[\alpha^{-1}\Gamma \alpha : \Gamma \cap 
\alpha^{-1}\Gamma \alpha] < \infty. $   

Let $\chi \in \rep(\Gamma,V).$ By definition, 
$\com(\Gamma) \subset \Gamma.$ However, in 
order to define non-trivial Hecke operators, we will need to 
start with an element $\alpha \in \com(\Gamma)
$ that lies outside of $\Gamma.$ We will also 
need the unitary representation $\chi$ to 
extend from $\Gamma$ to the set $\Gamma \alpha 
\Gamma.$ More precisely:

\begin{as} \label{asAlpha}
Throughout this paper we assume: 
\begin{enumerate}
\item The element $\alpha \in \com(\Gamma) \setminus \Gamma;$
\item The unitary representation $\chi \in \rep
(\Gamma,V)$ has a single-valued extension from 
$\Gamma$ to the set $\Gamma \alpha \Gamma,$ 
and satisfies the following properties:
\begin{enumerate}
\item $\chi(g_{1} \alpha g_{2}) = \chi(g_{1}) \chi(\alpha) \chi(g_{2})$ for $g_{1},g_{2}\in \Gamma;$
\item $\chi(\alpha)$ is an invertible linear map on the inner product space $V.$ 
\end{enumerate}
\end{enumerate}
\end{as}

See \cite[Pages 463-4]{Hejhal1} for the 
analogous Fuchsian case.

Let $d = 
[\Gamma : \Gamma \cap \alpha^{-1}\Gamma 
\alpha].$ It follows, by conjugation, that $[\alpha^{-1}\Gamma 
\alpha : \Gamma \cap \alpha^{-1}\Gamma \alpha] 
= d.$ 

Let the symbol $\bigsqcup$ denote disjoint union. Then 
$$\Gamma = \bigsqcup_{i=1}^{d} \left( \Gamma \cap \alpha^{-1}\Gamma \alpha \right) \epsilon_{i} \quad \text{iff} \quad \Gamma \alpha \Gamma = \bigsqcup_{i=1}^{d} \Gamma \alpha \epsilon_{i}.$$  

Setting $\alpha_{i} = \alpha \epsilon_{i}$ ($i=1 \dots d$), we define the operator $\M:\hs \mapsto \hs$ by
$$(\M f)(P) = \sum_{i=1}^{d} \chi(\alpha_{i})^{*}f(\alpha_{i} P).$$ 
The operator $\M$ does not depend on the above choice of coset representatives, and it is a bounded linear operator on $\hs$ (\cite[p. 467]{Hejhal1} or \cite{Strom1}).

The adjoint of $\M$ can be constructed as follows: Since $\alpha \in \com(\Gamma),$ and $\com(\Gamma)$ is a group; $\alpha^{-1}$ also satisfies Assumption~\ref{asAlpha}. As before, we can write  
$$\Gamma \alpha^{-1} \Gamma = \bigsqcup_{i=1}^{e} \Gamma \beta_{i}. $$ It follows that $e = d$ (\cite[p. 469]{Hejhal1}). Define $\M^{*}: \hs \mapsto \hs$ by
$$ (\M^{*}f)(z) = \sum_{i=1}^{d} \chi(\beta_{k}^{-1}) f(\beta_{k} P). $$
For every $f,g\in \hs$ we have (\cite[p. 470]{Hejhal1})
$$\left< \M f, g \right> = \left< f, \M^{*} g \right>. $$ 

\subsection{Interplay between $\lp$ and $\M$} 
For $P = z+rj,~P'=z'+r'j \in \HH$ set 
\beq \label{eqPointPair} 
\delta(P,P') \df  \frac{|z-z'|^{2}+r^{2}+ r'^{2}}{2rr'}. 
\eeq
It follows that $\delta(P,P') = \cosh(d(P,P'))$, where $d$ denotes the hyperbolic distance  in $\HH.$  Next, for $k \in \scz([1,\infty))$ a Schwartz-class function, define the \emph{point-pair invariant} $K(P,Q) \df k(\delta(P,Q)).$ Note that for all $\gamma \in \pc$ we have $K(\gamma P,\gamma Q) = K(P,Q). $
Define the Poincar\'{e} series 
$$
K_\Gamma(P,Q) \df \sum_{T\in\Gamma} \chi(T)K(P,T Q). $$ 
The series above converges absolutely and uniformly on compact subsets of
$\HH \times \HH$, and is the kernel of a bounded operator $\K : \hs \mapsto \hs.$  The Selberg theory allows us to \emph{diagonalize} $\K$ with respect to a basis of eigenfunctions of $\lp.$ 

We first need the following spectral expansion:   
Each $f \in \hs$ has an expansion of the form 
\begin{equation}\label{E:spectral}
f(P) = \sum_{m \in \D} \left< f,e_m \right> e_m(P). 
\end{equation}
Here the sum converges in the Hilbert space $\hs, $ and $\D$ is an indexing set of the eigenfunctions $e_m$ of $\lp$ with corresponding eigenvalues $\lambda_m.$

Now, let $h$ be the Selberg--Harish-Chandra transform of $k;$ explicitly, for $\lambda =1-s^{2} \in \CC,$  
\beq \label{eqSHC}
h(\lambda) = h(1-s^2) \df \frac{\pi}{s}
\int_{1}^{\infty}k\left(\frac{1}{2}\left(t+\frac{1}{t}\right)\right)
(t^{s}-t^{-s})\left(t-\frac{1}{t}\right)\,\frac{dt}{t};
\eeq
and let $g$ be the Fourier transform of $h:$
\beq \label{eqG}
g(x) = \frac{1}{2\pi} \int_{\RR} h(1+t^2)e^{-itx}\,dt. 
\eeq

For $v,w \in V $ let $ v \otimes \overline{w}$ be a linear operator in $V$ defined by    
$ v \otimes \overline{w}(x) = <x,w>v.  $
An immediate application of the spectral expansion \eqref{E:spectral} and the Selberg--Harish-Chandra transform \eqref{eqSHC} gives us (see \cite{Friedman1, Friedman2} or \cite[Equation 6.4.10, page 278]{Elstrodt}):
\begin{lem}
Let  $k \in \scz $ and  $h:\CC\ra\CC$ be the Selberg--Harish-Chandra  Transform
of $k.$ Then with $K_\Gamma$ defined above, we have 
\begin{equation}
\label{E:kernel expansion}
K_\Gamma(P,Q)   = \sum_{m \in \D} h(\lambda_m)e_m(P) \otimes 
\overline{e_m(Q)}
\end{equation}
The sum converges absolutely and uniformly on compact subsets 
of $\HH \times \HH$.
\end{lem}

Since the Laplace operator $\lp$ commutes with the action of $\pc$ on $\HH,$ it follows that $\lp$ commutes with $\M.$ From Selberg's original ideas (\cite{Selberg1}), the operator $\K = h(\lp),$ and thus $\M$ must also commute with $\K.$ A direct proof of this fact is given  in \cite[p. 473]{Hejhal1}:
\begin{lem} \label{lemCommu}
Let $\K$ be the linear operator defined from the kernel function $K_{\Gamma}(P,Q).$ Then, on $\hs$  
$$\K \M = \M \K, $$
and on $\hs \cap C^{2}(\HH)$
$$\M \lp = \lp \M. $$
\end{lem}  

For each eigenvalue $\lambda = \lambda_{m} ~ (m \in \D),$ let $A_{\lambda} \df A(\lambda, \Gamma, \chi)$ denote  the subspace of $\hs$ spanned by $\{ \lambda_{m}~|~ \lambda_{m} = \lambda   \}.$ Since $\lp$ commutes with $\M,$ it follows that the subspaces $A_{\lambda}$ are invariant under $\M.$ Let $m_{1}\dots m_{k_{\lambda}} \in \D $ generate $A_{\lambda}$ ($\lambda$ has multiplicity $k_{\lambda}$). Since $A_{\lambda}$ is invariant under $\M,$ the eigenvalues of $\left. \M \right|_{A_{\lambda}}$ can be listed as $\omega_{m_{1}}\dots \omega_{m_{k_{\lambda}}}.$ Next, we define   
\beq \label{eqTraceMod}
\omega(\lambda) \df \tr  \left. \M \right|_{A_{\lambda}} = \sum_{j=1}^{k_{\lambda}} \omega_{m_{j}}. 
\eeq

Decomposing $\hs$ into invariant subspaces $A_{\lambda},$ for each $m \in \D$ we pair $\lambda_{m}$ and $\omega_{m}.$ Since $\M$ is a bounded operator, we conclude that the $\omega_{m}$ $(m \in \D)$ are uniformly bounded.

\section{The Selberg trace formula for Hecke operators} \label{secSelTra}
The standard Selberg trace formula is an explicit evaluation of both the integral and spectral of the operator $\K.$ In our case of interest we evaluate the trace of $\K \M.$
\subsection{The kernel of operator $\K \M$}
The first step is to construct a Poincar\'{e} series for the kernel of $\K \M.$ 

For all $f \in \hs,$ by Lemma~\ref{lemCommu}, and by definition, 
$$\K \M f(P) = \M \K f(P) = \int_{\F} \M_{P} K_{\Gamma}(P,Q) f(Q) \, dv(Q). $$
Here $\M_{P} K_{\Gamma}(P,Q) $ is the action of $\M,$ treating $K_{\Gamma}(P,Q)$ as a function of $P.$ Set  
\beq
K_{\M}(P,Q) \df \sum_{T \in \Gamma \alpha^{-1} \Gamma } \chi(T^{-1})^{*} K(P,T Q),  
\eeq  
it\footnote{Note that $\chi$ is not necessarily unitary when extended to $\Gamma \alpha \Gamma,$ so $\chi(T^{-1})^{*}$ can not be simplified. } then follows (\cite[p. 476]{Hejhal1}) that 
$K_{\M}(P,Q) = \M_{P} K_{\Gamma}(P,Q). $ Hence:
\begin{lem}
For all $f \in \hs,$ 
$$\K \M f(P) = \M \K f(P) =  \int_{\F} K_{\M}(P,Q) f(Q) \, dv(Q). $$
\end{lem}

Upon applying $\M$ to \eqref{E:kernel expansion}, we obtain
\beq
K_{\M}(P,Q) = \sum_{m \in \D} h(\lambda_m) \left(\M e_m(P)\right) \otimes 
\overline{e_m(Q)}.
\eeq 
The Selberg trace formula will result from a careful evaluation of 
$$ \tr(\K \M) \df \int_{\F} \tr_{V} K_{\M}(P,P)\, dv(P). $$

\subsection{Spectral trace}
Since $K_{\M}(P,Q)$ is a continuous function on the compact set $\F \times \F,$ the integral 
$$ \int_{\F} \tr_{V} K_{\M}(P,P)\, dv(P) $$ converges absolutely. Since $\M$ is bounded, the derivation of the standard Selberg trace formula for cofinite Kleinian groups \cite[\S6.5]{Elstrodt} implies that
\beq \label{eqSpecTra}
\int_{\F} \tr_{V} K_{\M}(P,P)\, dv(P) =  \sum_{m \in \D} h(\lambda_m) \omega_{m}.
\eeq

\subsection{Integral trace}
Selberg's clever method of evaluating the integral \eqref{eqSpecTra} entails decomposing the Poincar\'{e} series $K_{\M}(P,Q)$ into conjugacy classes of various types of isometries (the types being: elliptic, parabolic, hyperbolic, \dots). A similar method works here, though conjugacy classes must be defined in the context of the set $\Gamma \alpha^{-1} \Gamma.$ One also must deal with centralizer subgroups of various isometries, which lead to complications not found in the standard Selberg trace formula derivations.

For $T \in \Gamma \alpha^{-1} \Gamma$ set
$$\{ T \} \df \{ T \}_{\Gamma} \df \{ \sigma^{-1} T \sigma~|~\sigma \in \Gamma  \},$$
and set
$$\C(T) \df \C_{\Gamma}(T) \df \{ \sigma \in \Gamma~|~ \sigma T = T \sigma  \}. $$

It now follows (\cite[pp. 480--482]{Hejhal1}) that 
\beq \label{eqSelTrick}
\tr(\K \M) = \sum_{\{ T\} } \tr_{V}\!\left( \chi(T^{-1})^{*} \right) \int_{\F(\C(T))}K(P, T P)\, dv(P).
\eeq
Here the sum is over all conjugacy classes of $\Gamma \alpha^{-1} \Gamma$ with respect to $\Gamma;$ $\F(\C(T))$ is a fundamental domain for the subgroup $\C(T).$ 

Since the integral in \eqref{eqSelTrick} is invariant when $T$ is conjugated within $\pc,$ it suffices to choose a conjugate of $T$ which has the simplest form to allow one to explicitly evaluate the following integral:
For $T \in \Gamma \alpha^{-1} \Gamma$ set 
\beq \label{eqIntSel}
I(T) \df I_{\Gamma}(T) = \int_{\F(\C(T))}K(P, T P)\, dv(P).
\eeq
There are four cases to consider:
\begin{enumerate}
\item $T = \text{id}$ (this case does not occur by Assumption~\ref{asAlpha})
\item $T$ is loxodromic
\item $T$ is elliptic
\item $T$ is parabolic.
\end{enumerate}
In fact, we must first compute the structure of each centralizer subgroup $\C(T).$

See \mbox{\cite[chap. 1]{Elstrodt}} for details on the action of $\pc$ on  $\HH,$ and for the definitions of loxodromic, parabolic, \dots

\subsection{Loxodromic case} Suppose that $T \in \Gamma \alpha^{-1} \Gamma$ is loxodromic. We can then assume (replacing $T$ by its unique conjugate within $\pc$) that $T$ has the form  
$$
T=
\left(\begin{array}{cc}
a(T) & 0\\
0 & a(T)^{-1}
\end{array}\right) $$
such that $a(T)\in\CC$ has $|a(T)|>1$.  Let $N(T)$  denote  the \emph{norm} of $T,$   defined by  $$N(T) \df |a(T)|^{2}.$$
For $z + rj \in \HH$ 
\beq \label{eqLoxAct}
T(z+rj) = a(T)^{2}z + N(T) rj.  
\eeq

The element $T$ fixes two points: $0$ and $j \infty \in \widehat{\CC}.$ It follows from elementary calculations that any $\sigma \in \C(T) \subset \Gamma$ must be of the form 
$$\left(\begin{array}{cc}
u & 0\\
0 & u^{-1}
\end{array}\right). 
$$ 
Since $\Gamma$ is discrete and cocompact, it follows (\cite[chap. 5]{Elstrodt}) that $\C(T)$ has at least one element of minimal norm. Choose one, say $T_{0}.$ The other minimal norm elements are obtained by multiplying by the elliptic elements of the finite cyclic elliptic subgroup  $\en$ of order $m(T), $ generated by an element $E_{T_0}. $   We now have  
$$
\C(T) = \langle T_{0} \rangle \times \en. 
$$
Here $\langle  T_{0} \rangle = \{\, T_{0}^{n} ~ | ~ n \in\ZZ ~ \}.$ The elliptic element $ E_{T_0}$ is conjugate in $\pc$ to an element of the form 
$$\left(\begin{array}{cc}
\zeta(T_0) & 0 \\
0 & \zeta(T_0)^{-1}
\end{array}\right), $$ 
where here $\zeta(T_0)$ is a primitive $2m(T)$-th root of unity.

Note that (by definition) if $T \notin \Gamma,$ $T \notin \C(T).$ It seems conceivable then that $\C(T)$ could be a finite group. However, in the next lemma we will show that $N(T_{0}) > 1$ which implies that $\C(T)$ contains an infinite cyclic loxodromic subgroup\footnote{This is closely analogous to the Fuchsian group case, where the centralizer of a hyperbolic element is infinite \cite[p. 483]{Hejhal1}. Another analogous result: ``\dots In dimension 3, the centralizer of an elliptic element of a group with compact quotient never is a finite cyclic group. '' \cite[p. 194]{Elstrodt}.Meanwhile in dimension two, it is a finite cyclic group.   }.   

\begin{lem} \label{lemMoreThanEllip}
Let $T \in \Gamma \alpha^{-1} \Gamma.$ Then $\C(T)$ contains a loxodromic element $T_{0}$ of minimal norm $N(T_{0}) > 1$. 
\end{lem}
\begin{proof}
Let $k \in \scz[1,\infty)$ be a smooth, positive bump function with $k \df 1$ on $[1,2],$ (zero everywhere else) and let $\chi$ be the trivial representation. The right-hand side of \eqref{eqSpecTra} converges absolutely, and by positivity the left side does, too. Hence \eqref{eqSelTrick} and \eqref{eqIntSel} imply that 
\beq \label{eqInter82382} 
\int_{\F(\C(T))}K(P, T P)\, dv(P) < \infty. 
\eeq 
For $\varphi \in \RR$ set
$$R(\varphi) =  \left(\begin{array}{cc}
e^{\frac{i \varphi}{2}} & 0 \\
0 & e^{-\frac{i \varphi}{2}}
\end{array}\right).  
$$
$R(\varphi)$ is a rotation on $\HH$ that fixes $0,\infty j \in \widehat{\CC}.$ For $z+rj \in \HH,$ 
$$
R(\varphi)(z+rj) = e^{i\varphi}z + rj.
$$
Next, suppose there are not loxodromic elements in $C(T),$ that is, that  $\C(T) = \langle R(\varphi)  \rangle,$ where $\varphi \df \frac{2\pi}{m},$ and $m = |R(\varphi)| = |\C(T)|. $
A fundamental domain for $\C(T)$ is given by (see \cite[p. 193]{Elstrodt})
$$\F_{1} = \left\{  \rho e^{i\varphi} + rj ~|~ r > 0,~\rho \geq 0,~ 0 \leq \varphi \leq \frac{2\pi}{m}   \right\}. $$
We claim that 
\beq \label{eqInter3729} 
\int_{\F_{1}}K(P, T P)\, dv(P) = \infty. 
\eeq 
To see this, first note that using \eqref{eqLoxAct} and \eqref{eqPointPair}, we can simplify  $K(P, T P).$ We obtain
$$ 
\int_{\F_{1}}K(P, T P)\, dv(P) = 
\int_{\F_{1}}k\left( \frac{|a(t)^{2}-1|^{2}|z|^{2} + (N(T)^{2}+1)^{2}r^{2}}{2N(T) r^{2}}   \right)   \, \frac{dx\, dy\, dr}{r^{3}}.
$$
Applying the elementary substitution $x \mapsto rx, ~ y \mapsto ry;$ the integral becomes
\begin{multline*}
\int_{\F_{1}}k\left( \frac{|a(t)^{2}-1|^{2}|z|^{2} + (N(T)^{2}+1)^{2}}{2N(T) }   \right)   \, \frac{dx\, dy\, dr}{r}  \\ = \int_{0}^{\infty} \frac{dr}{r} \int_{z(\F_{1})} k\left( \frac{|a(t)^{2}-1|^{2}|z|^{2} + (N(T)^{2}+1)^{2}}{2N(T) }   \right)   \, dx\, dy = \infty.  
\end{multline*}
Here $z(\F_{1})$ is the standard projection of $\F_{1} \subset \HH$ onto $\CC.$  
Hence $\C(T)$ must contain at least one loxodromic element, and since $\Gamma$ is discrete, $\C(T)$ must contain a minimally normed element.
\end{proof}
An element $R$ in $\Gamma$ is called \emph{primitive} if there \emph{does not} exist $S \in \Gamma$ with $S^{m} = R$ for some positive integer $m.$ Since the $T_{0}$ from Lemma~\ref{lemMoreThanEllip} has minimal norm in $\C(T),$ it is necessarily primitive.  

Now that we see that the structure of $\C(T)$ is identical to the case where $T \in \Gamma$ (\cite[p. 193]{Elstrodt}), hence we can evaluate the integral \eqref{eqIntSel}, obtaining:

\begin{lem} \label{lemLoxo}
Let \emph{$T \in \com(\Gamma)$} be loxodromic (hyperbolic elements are considered loxodromic). Let 
$$
\C(T) = \langle T_{0} \rangle \times \en 
$$
with $T_{0}$ loxodromic\footnote{Note that there is \emph{no} clear relationship between $T$ and $T_{0}$ other than the fact that they have the same fixed points and commute. It seems reasonable to conjecture that $T=T_{0}^{k/n}$ for $k,n \in \NN.$ } and primitive, and $\en$ a finite cyclic subgroup of elliptic elements. Then 
$$
\int_{\F(\C(T))}K(P, T P)\, dv(P) =  \frac{g( \log N(T))
\log N(T_0)}{|\E(T)| \, |a(T) - a(T)^{-1}|^2}.
$$
Here, $g$ comes from \eqref{eqG} and \eqref{eqSHC}.
\end{lem} 

\subsection{Elliptic case} \label{secEllCas}
Suppose that $R\in \Gamma \alpha^{-1} \Gamma$ is elliptic. We can then assume (replacing $R$ by its unique conjugate within $\pc$) that $R$ has the form 
$$R= R(\varphi) =  \left(\begin{array}{cc}
e^{\frac{i \varphi}{2}} & 0 \\
0 & e^{-\frac{i \varphi}{2}}
\end{array}\right)  \quad (\varphi \in \RR).  
$$
$R(\varphi)$ is a rotation on $\HH$ that fixes $0,\infty j \in \widehat{\CC}.$ For $z+rj \in \HH,$ 
$$
R(\varphi)(z+rj) = e^{i\varphi}z + rj.
$$
An elementary calculation shows that 
$$
\left(\begin{array}{cc}
a & b\\
c & d
\end{array}\right) \in \pc 
$$
commutes with $R(\varphi)$ iff one of the two conditions are satisfied: 
\begin{enumerate}
\item $b=c=0$
\item $\exp{\frac{i\varphi}{2}} = \pm i $ and $a=d=0.$  
\end{enumerate}
In both cases, a similar argument to Lemma~\ref{lemMoreThanEllip} (see \cite[pp. 193--194]{Elstrodt}) shows that $\C(R)$ contains loxodromic elements. In particular $\C(R)$ contains a loxodromic element of minimal norm $T_{0}.$ Following \cite[pp. 193--197]{Elstrodt} and generalizing to $R\in \com(\Gamma)$ we  have (\cite[Thm. 5.2.1]{Elstrodt})\footnote{We modified their hypothesis to allow $R$ to be in $\com(\Gamma)$ instead of in $\Gamma.$ We are justified by the \emph{proof} of Lemma~\ref{lemMoreThanEllip}.}

\begin{lem} \label{lemEllStr}
Let \emph{$R\in \com(\Gamma)$} be elliptic. Then there exists a primitive elliptic element $R_{0} \in \C(R) \subset \Gamma $ (it is possible that $R_{0} = \emph{id}$), an element $T_{0} \in \C(R)$ loxodromic, of minimal norm in $\C(R)$ so that one of the following cases holds:
\begin{enumerate}
\item Either all the elliptic elements of $\C(R)$ are contained in  $\E(R) \df \left<  R_{0} \right>,$ and 
$$\C(R) = \left< T_{0}\right> \times \E(R).$$
\item Or $R$ is elliptic of order 2, and there exists an elliptic element $S \in \C(R)$ so that 
$$\E(R) = \left< R_{0} \right> \cup \left< R_{0} \right>{S}$$ and    
$$\C(R) = \left< T_{0}\right> \times \E(R). $$ 
\end{enumerate}
\end{lem}
See \cite[pp. 191-198]{Elstrodt} for more details.

Now that we see that the structure of $\C(T)$ is identical to the case where $T \in \Gamma$ (\cite[p. 193]{Elstrodt}), we can evaluate the integral \eqref{eqIntSel}, obtaining:

\begin{lem} \label{lemEllip}
Let \emph{$R \in \com(\Gamma)$} be elliptic. Let 
$$
\C(R) = \langle T_{0} \rangle \times \E(R) 
$$
as in Lemma~\ref{lemEllStr}. Then 
$$
\int_{\F(\C(R))}K(P, R P)\, dv(P) =  \frac{  g(0) \log N(T_0) }{  |\E(R)| \, |(\tr(R))^{2}-4|}.
$$
Here, $g$ comes from \eqref{eqG} and \eqref{eqSHC}.
\end{lem}

We will also need the following:
\begin{lem} \label{lemFinEll}
Let \emph{$R \in \com(\Gamma)$} be elliptic. Then there are only finitely many elliptic conjugacy classes of $\Gamma \alpha^{-1} \Gamma,$ where conjugacy is defined with respect to $\Gamma.$
\end{lem}
\begin{proof}
Let $k \in \scz[1,\infty)$ be non-negative and let $\chi$ be the trivial representation. By \eqref{eqSelTrick} and \eqref{eqSpecTra}, 
$$ \sum_{ \{ R \}_{\text{ell}}  } \int_{\F(\C(R))}K(P, R P)\, dv(P) < \infty. $$
Hence by Lemma~\ref{lemEllip}, 
$$ g(0) \sum_{ \{ R \}_{\text{ell}}  } \frac{\log N(T_0) }{  |\E(R)| \, |(\tr(R))^{2}-4|} < \infty. $$
Now, since $N(T_{0}) \geq \delta > 1$ (independent of $T_{0}$) and since $-2 < \tr(R) < 2,$ in order for there to be infinitely many conjugacy classes, it is necessary for $|\E(R)|$ to be unbounded as $R$ goes through all conjugacy classes. But by Lemma~\ref{lemEllStr} that would imply the order of each $R_{0}$ is unbounded; a contradiction since cofinite  Kleinian groups have only finitely many conjugacy classes of elliptic elements (and two elements in the same conjugacy class have the same order; see\cite{Elstrodt}).          
\end{proof}

\subsection{Parabolic case} \label{secPara} Suppose that $T \in \com(\Gamma)$ is parabolic. We will soon see that this case can not occur when $\Gamma$ is assumed to be cocompact. That is: $\com(\Gamma)$ \emph{contains no parabolic elements when $\Gamma$ is cocompact.}
To see this, suppose that $T$ is already conjugated into the simple form 
$$
T = \left(\begin{array}{cc}
1 & \tau \\
0 & 1
\end{array}\right). 
$$
Then, an elementary calculation shows that 
$$
S \df \left(\begin{array}{cc}
a & b\\
c & d
\end{array}\right) \in \pc 
$$
commutes with $T$ iff $c=0$ and $a=d = \pm 1.$ By definition, $S$ must be parabolic (or the identity). But, since $\Gamma$ is cocompact, it contains no parabolic elements. So $\C(T) = \text{id}.$  An argument similar to Lemma~\ref{lemMoreThanEllip} shows that this can not happen. Hence there are no parabolic elements in $\com(\Gamma) \setminus \Gamma.$

%that for  $k \in \scz[1,\infty)$ a smooth, positive bump function with $k \df 1$ on $[1,2]$ (zero everywhere else),  $\chi$ the trivial representation. Once again, the right-hand-side of \eqref{eqSpecTra} converges absolutely, so does the left---and in particular, by \eqref{eqSelTrick} and \eqref{eqIntSel}, 
%$$ 
%\int_{\F(\C(T))}K(P, T P)\, dv(P) < \infty. 
%$$ 
%Since $\C(T) = \text{id};$ $\F(\C(T)) = \HH.$ An elementary calculation (like in Lemma~\ref{lemMoreThanEllip}) shows that the above integral diverges. Thus, we finally conclude that $T$ can not be parabolic.

\subsection{The Selberg trace formula}
By putting together \eqref{eqSelTrick}, Lemma~\ref{lemLoxo}, Lemma~\ref{lemEllip}, Assumption~\ref{asAlpha}, and the results of \S\ref{secPara}; and by applying standard approximation techniques in order to allow more general growth conditions on the function $k$ (and hence on $h$ and $g$) (see \cite[pp. 32--34]{Hejhal1} for the details), we obtain:
    
\begin{thm}{(Selberg trace formula)} \label{T:Selberg}
Let $\Gamma $ be a cocompact Kleinian group. Let \emph{$\alpha \in\com(\Gamma) \setminus \Gamma,$} and let  $\chi$ be a representation satisfying Assumption~\ref{asAlpha}. Let  $h$ be a holomorphic function on  $ \{ s \in \CC \, | \, |\I(s)| < 2+ \delta \}$ for some $\delta > 0,$ satisfying $ h(1+z^2) = O( 1+|z|^2)^{3/2 - \epsilon}) $ as  $|z| \ra \infty,$ and let
$$ g(x) = \frac{1}{2\pi} \int_{\RR} h(1+t^2)e^{-itx}\,dt. $$  Then
\begin{multline*}
\sum_{m \in \D} h(\lambda_m) \omega_{m} = \sum_{ \{R \}_{\text{\emph{ell}}}} \frac{\tr_{V}\!\left( \chi(R^{-1})^{*} \right)   g(0) \log N(T_0) }{  |\E(R)| \, |(\tr(R))^{2}-4|}  + 
\sum_{\{ T \}_\text{\emph{lox}} } \frac{\tr_{V}\!\left( \chi(T^{-1})^{*} \right) g(\log N(T))}{\oen |a(T)-a(T)^{-1}|^{2}}\log N(T_{0})
\end{multline*}
\end{thm}
Here, $ \{ \lambda _{m} \}_ {m \in \D } $ are the eigenvalues of $\lp$ counted with multiplicity, and $\omega_{m}$ are the eigenvalues of $\M$ (the Hecke operator associated to $\alpha$) with the convention of Equation~\ref{eqTraceMod}.  The summation with respect to $\{R \}_\text{ell} $ extends over the finitely many $\Gamma-$conjugacy classes of  elliptic elements  $R \in \Gamma \alpha^{-1} \Gamma,$ and for such a class, $N(T_0)$ is the minimal norm of a hyperbolic or loxodromic element of the centralizer $\mathcal{C}(R) \subset \Gamma;$ $\ren$ is the maximal finite subgroup contained in $ \mathcal{C}(R)$ (see Lemma~\ref{lemEllip} for more details).  The summation with respect to $\{ T \}_\text{lox}$ extends over the  $\Gamma-$conjugacy classes of hyperbolic or loxodromic elements of $\Gamma \alpha^{-1} \Gamma,$ $T_0$ denotes a primitive hyperbolic or loxodromic element of minimal norm in\footnote{Note that there is no clear relationship between $T \in \Gamma \alpha^{-1} \Gamma$ and $T_{0} \in \Gamma.$} $\mathcal{C}(T) \subset \Gamma;$  the element $T$ is conjugate in $\pc$ to the transformation  described by the diagonal matrix with diagonal entries $a(T), a(T)^{-1}$ with $|a(T)| > 1, $ and\footnote{Please note that there is a typographical error in the loxodromic and non cuspidal elliptic terms in \cite{Elstrodt} Theorem 6.5.1; both terms are missing a factor of $\frac{1}{4 \pi}. $} $N(T) =   |a(T)|^2; $ $\en$ is the finite cyclic elliptic subgroup of $\mathcal{C}(T)$  (see Lemma~\ref{lemLoxo} for more details).  The sum over elliptic elements is finite, and all other sums converge absolutely.

Please keep in mind that the conjugacy classes $\{R\}_{\text{ell}}$ and $\{T\}_{\text{lox}},$  the operator $\M,$ and the eigenvalues $\omega_{m}$ all depend on $\alpha \in \com(\Gamma) \setminus \Gamma.$ 

\subsection{Self-Adjoint Hecke Operators}
Throughout this work we were under Assumption~\ref{asAlpha}. By imposing stronger conditions on $\alpha$ and $\chi$ we can ensure that $\M$ is self-adjoint.

\begin{as} \label{asSelf}
In addition to Assumption~\ref{asAlpha}
\begin{enumerate}
\item $\Gamma \alpha \Gamma = \Gamma \alpha^{-1} \Gamma$
\item $\chi(\alpha^{-1}) = \chi(\alpha)^{*}. $
\end{enumerate}
\end{as}   

Assumption~\ref{asSelf} implies that $\M$ is self-adjoint. Since $\M$ commutes with $\lp$ we can choose a simultaneous diagonalization $\{ e_{m} \}_{m \in \D}$ so that
$$ \lp e_{m} = \lambda_{m} e_{m} \quad \text{and} \quad \M e_{m} = \omega_{m} e_{m}.$$ 

\section{Application of the Trace Formula to the distribution $\omega_{m}.$}
In the Selberg trace formula (Theorem~\ref{T:Selberg}), the eigenvalues $\lambda_{m}$ of $\lp$ \emph{do not} depend on $\M.$ The distribution of these eigenvalues has been well studied \cite{Huber1, Selberg1, Hejhal1, Elstrodt, Venkov}. More specifically, in the case of cocompact Kleinian groups we have 
\begin{thm}\cite[p. 308]{Elstrodt}
Let $\Gamma$ be a cocompact Kleinian group. Then
$$\sum_{m \in \D} e^{-\lambda_{m}t} \sim \frac{\vol(\Gamma)}{8 \pi^{3/2}}t^{-3/2} \quad \text{as $t \ra 0^{+}.$} $$
Here $\vol(\Gamma)$ is the volume of a fundamental domain for the action of $\Gamma$ in $\HH.$
\end{thm}

For  $\Gamma$  a cocompact Kleinian group with $\alpha \in \com(\Gamma) \setminus \Gamma,$ define $E_{\alpha}$ the  \emph{elliptic number of $\Gamma$ with respect to $\alpha$} by
\beq \label{eqEllNum}
E_{\Gamma}^{\alpha} \df \sum_{ \{R \}_{\text{{ell}}}} \frac{ \log N(T_0) }{  |\E(R)| \, |(\tr(R))^{2}-4|} 
\eeq
(see \S\ref{secEllCas} for the definitions of the above notation).

Our result is as follows:
\begin{thm} \label{thmHeckLap}
Let $\Gamma$ be a cocompact Kleinian group with $\alpha \in \com(\Gamma) \setminus \Gamma.$  Then
$$  \sum_{m \in \D} \omega_{m} e^{-\lambda_{m}t} = \frac{E_{\Gamma}^{\alpha}}{\sqrt{4\pi t}} + O(\sqrt{t}) \quad \text{as $t \ra 0^{+}.$} $$
Here $E_{\Gamma}^{\alpha}$ is the elliptic number of $\Gamma$ with respect to $\alpha,$ $\{ \lambda_{m} \}_{m \in \D } $ are the eigenvalues of $\lp$ counted with multiplicity, and $\omega_{m}$ are the eigenvalues of $\M$ (the Hecke operator associated to $\alpha$) with the convention of Equation~\ref{eqTraceMod}.
\end{thm}
\begin{proof} 
It follows from Lemma~\ref{lemSmallest} (below) that the loxodromic sum in \eqref{eqHeat} is bounded by $O\left(t^{-1/2}\exp(-c/t)\right)$ as $t \ra 0^{+}.$ Hence 
\begin{multline*}
\sum _{m \in \D } \omega_{m} e^{-t\lambda _{m}} = \\ \frac{\exp(-t)}{\sqrt{4\pi t}} \sum_{ \{R \}_{\text{{ell}}}} \frac{ \log N(T_0) }{  |\E(R)| \, |(\tr(R))^{2}-4|}  + O\left(t^{-1/2}\exp(-c/t)\right)  = \frac{E_{\Gamma}^{\alpha}}{\sqrt{4\pi t}} + O(\sqrt{t}) \quad \text{as $t \ra 0^{+}.$}   \end{multline*}
\end{proof}

\begin{lem} \label{lemSmallest}
Let $\Gamma$ be a cocompact Kleinian group with $\alpha \in \com(\Gamma) \setminus \Gamma.$ Then there exists a constant $c_{0} > 1$ so that $N(T) > c_{0}$ for all conjugacy class $\{ T \} \in \{T\}_{\text{lox}}.$
\end{lem}
\begin{proof}
Recall that $\{T\}_{\text{lox}}$ is the set of conjugacy classes of loxodromic elements of 
$\Gamma \alpha^{-1} \Gamma,$ where conjugacy is taken with respect to the group $\Gamma.$ An application of the Selberg trace formula, for the pair of functions 
$$h(z)=\exp(-zt), \quad
g(r)=\frac{\exp(-t)}{\sqrt{4\pi t}}\exp\left(-\frac{r^2}{4t}\right),
$$ 
yields 
\begin{multline} \label{eqHeat}
\sum _{m \in \D } \omega_{m} e^{-t\lambda _{m}} = \\ \frac{\exp(-t)}{\sqrt{4\pi t}} \sum_{ \{R \}_{\text{{ell}}}} \frac{ \log N(T_0) }{  |\E(R)| \, |(\tr(R))^{2}-4|}  + \frac{\exp(-t)}{\sqrt{4\pi t}} \sum_{\{ T \}_\text{{lox}} } \frac{ \exp\left(-\frac{(\log N(T) )^2}{4t}\right)}{\oen |a(T)-a(T)^{-1}|^{2}}\log N(T_{0}).
\end{multline}
For $t > 0$ the left hand side of \eqref{eqHeat} converges absolutely. Fix $t = 1.$ Now, the sum over loxodromic terms is comprised of positive terms. Hence each term of the form   
$$ \frac{ \exp\left(-\frac{(\log N(T) )^2}{4}\right)}{\oen |a(T)-a(T)^{-1}|^{2}}\log N(T_{0}) $$ must be bounded by a constant \emph{independent of $T$ and $T_{0}.$} Since $\Gamma$ is a discrete group, there exists a constant $d_{0} > 1$ so that $N(T_{0}) > d_{0}$ for all loxodromic $T_{0} \in \Gamma.$  Hence $\log N(T_{0})$ is bounded above zero, uniformly (since $T_{0} \in \Gamma$). It thus follows that $|a(T)-a(T)^{-1}|$ must be uniformly bounded above zero. In other words, there is a constant $c_{0} > 1,$ independent of $T,$ with $|a(T)| > c_{0}.$ Noting that $N(T)=|a(T)|^{2},$ the Lemma follows.  
\end{proof}

Theorem~\ref{thmHeckLap} immediately implies:
\begin{thm} 
Let $\Gamma$ be a cocompact Kleinian group with $\alpha \in \com(\Gamma) \setminus \Gamma.$ Suppose that the set $\Gamma \alpha^{-1} \Gamma$ contains no elliptic elements. Then
$$  \sum_{m \in \D} \omega_{m} e^{-\lambda_{m}t} = 
 O\left(t^{-1/2}\exp(-c/t)\right)  
 \quad \text{as $t \ra 0^{+}.$} $$
Hence
$$
\lim_{t \ra 0^{+}} \sum_{m \in \D} \omega_{m} e^{-\lambda_{m}t}   = 0.$$
Here $\{ \lambda_{m} \}_{m \in \D } $ are the eigenvalues of $\lp$ counted with multiplicity, and $\omega_{m}$ are the eigenvalues of $\M$ (the Hecke operator associated to $\alpha$) with the convention of Equation~\ref{eqTraceMod}.
\end{thm}

\section{Analogues of Huber's Theorem for Hecke Operators}
Following \cite[p. 202]{Elstrodt}, we define the length spectrum of $\Gamma \alpha^{-1} \Gamma.$ For loxodromic $T_{j} \in \Gamma \alpha^{-1} \Gamma$ set $\mu_{j} = \log N(T_{j}).$  The \emph{length spectrum} of $\Gamma \alpha^{-1} \Gamma$ is defined to be (see \S\ref{secSelTra} for the notation)    
$$ \mathcal{L}_{\Gamma}^{\alpha} \df \left( \mu_{j}, \sum_{ \substack{ \{ T \}_\text{\emph{lox}} \\ \log N(T)=\mu_{j} } }  \frac{\log N(T_{0}) }{\oen |a(T)-a(T)^{-1}|^{2}}  \right)_{j \geq 1.} $$ 

In the two-dimensional case the length spectrum is simply comprises the lengths of closed geodesics. Here we really need the \emph{complex lengths} $a(T)$ and the order of the elliptic, finite subgroup of $\C(T),$ $\oen.$  

We define the eigenvalue spectra of $\lp$ and $\M$ ($\M$ is defined from $\alpha \in  \com(\Gamma) \setminus \Gamma$) by
$$\mathcal{S}_{\Gamma}^{\alpha} = \left(\lambda_{j}, \omega(\lambda_{j}) \right)_{j \in \D^{*}};$$
here the symbol $D^{*}$ means that we do not count with multiplicity: recall that $\omega(\lambda)$ is the trace of $\M$ on the invariant subspace generated by all eigenfunctions (of $\lp$) with eigenvalue $\lambda,$ so multiplicity is already encoded into $\omega(\lambda).$

\begin{thm} Let $\Gamma, \Gamma'$ be cocompact Kleinian groups with $\alpha \in \com(\Gamma) \setminus \Gamma,~\alpha' \in \com(\Gamma') \setminus \Gamma'. $ Then the following hold:

\begin{enumerate}
\item Suppose that $\mathcal{S}_{\Gamma}^{\alpha}$ and $\mathcal{S}_{\Gamma'}^{\alpha'}$ agree up to at most finitely many terms. Then $$E_{\Gamma}^{\alpha} = E_{\Gamma'}^{\alpha'}$$ 
$$\mathcal{S}_{\Gamma}^{\alpha} = \mathcal{S}_{\Gamma'}^{\alpha'}$$  
$$\mathcal{L}_{\Gamma}^{\alpha} = \mathcal{L}_{\Gamma'}^{\alpha'}.$$
\item Suppose that $\mathcal{L}_{\Gamma}^{\alpha}$ and $\mathcal{L}_{\Gamma'}^{\alpha'}$ agree up to at most finitely many terms. Then $$E_{\Gamma}^{\alpha} = E_{\Gamma'}^{\alpha'}$$ 
$$\mathcal{S}_{\Gamma}^{\alpha} = \mathcal{S}_{\Gamma'}^{\alpha'}$$  
$$\mathcal{L}_{\Gamma}^{\alpha} = \mathcal{L}_{\Gamma'}^{\alpha'}.$$
\end{enumerate}
\end{thm}
\begin{proof}
Upon applying the Selberg trace formula (with trivial representations) to the pair of functions,  
$$h(w)=\frac{1}{s^2+w-1} - \frac{1}{B^2+w-1} $$ 
$$ g(x) = \frac{1}{2s}e^{-s|x|} - \frac{1}{2B}e^{-B|x|}, $$ where $1 < \R(s) < \R(B),$  we obtain
\begin{multline} \label{eqLogDer}
\frac{1}{2s} \sum_{ \{ T \}_{\text{lox}}}  \frac{ \log N(T_{0})}{m(T)|a(T)-a(T)^{-1}|^{2}}N(T)^{-s} 
 -\frac{1}{2B}\sum_{ \{ T \}_{\text{lox}}} \frac{  \log N(T_{0})}{m(T)|a(T)-a(T)^{-1}|^{2}}N(T)^{-B} \\
= \sum_{n \in D^{*}} \omega(\lambda_{n}) \left(\frac{1}{s^2 - s_n^2} - \frac{1}{B^2 - s_n^2}  \right) 
 - \left( \frac{1}{2s}-\frac{1}{2B}\right) \sum_{ \{R \}_\text{ell}}\frac{\log N(T_{0})}{  |\E(R)| \, |(\tr(R))^{2}-4|}.
\end{multline}
Here $s_{n}^{2} = 1-\lambda_{n}$ and $s_{n}$ is chosen to lie in $ \{ z \in \CC~|~\I{z}\geq 0, \R{z} \geq 0~\}.$ 
The Theorem follows by comparing asymptotics. Simply apply the proof of the \cite[Theorem 3.3, p. 203]{Elstrodt} replacing \cite[Equation 5.2.35, p. 198]{Elstrodt} with our \eqref{eqLogDer}. Note that our symbol $w(\lambda_{n})$ is the analogue of the multiplicity of $\lambda_{n}$ as considered in \cite{Elstrodt}.  
\end{proof}

By applying the above Theorem to the case of a \emph{single} fixed group with two different elements $\alpha, \alpha' \in \com(\Gamma) \setminus \Gamma, $ we obtain:

\begin{cor}
Let $\Gamma$ be a cocompact Kleinian groups with $\alpha, \alpha' \in \com(\Gamma) \setminus \Gamma. $ Let $\omega(\lambda_{m}) = \omega'(\lambda_{m})$ for all but at most finitely many $m \in \D^{*}.$ Then $\omega(\lambda_{m}) = \omega'(\lambda_{m})$ for all $m \in \D^{*}.$ 
\end{cor}

%%%%%%%%%%%%%%%%%%%%%%
\bibliographystyle{amsalpha}
\bibliography{mod}

\end{document}